
\input amstex

\magnification =1200
\documentstyle{amsppt}
\refstyle{A}
\widestnumber\key{ACGH22}

\hsize =16.0truecm
\baselineskip =15.8truept
 \vsize =23.58truecm
\parskip 4pt
\hcorrection{0.17truein}


\NoBlackBoxes
\NoRunningHeads

\topmatter
 \title  Complex Contact threefolds and their Contact Curves \endtitle
 \author  Yun-Gang Ye \\  \\  Department of Mathematics \\
 Texas A\&M University \\ College Station, TX 77843 \endauthor
\email ygye@math.tamu.edu \endemail
\endtopmatter

\TagsAsMath

\def\E{\Cal E}

\def\h#1#2{H^{#1}\left(#2\right)}
\def\cmt{\Omega_{X/M}^1}

\def\cpc{{\Bbb C}{\Bbb P}^3}
\def\cpb{{\Bbb C}{\Bbb P}^2}

\def\ppb{{\Bbb P}\left(\Omega^1_{\cpb}\right)}
 \def\ct{\Omega_X^1}
\def\pt{\Omega_{\cpc}^1}
\def\cpbd{{\Bbb C}{\Bbb P}^{2*}}
\def\pp{\cpb\times \cpbd}
\def\li{\ell_{\infty}}
\def\lid{\li^{*}}
\def\op_{{\Cal O}}
\def\hb{\Cal H}
\def\hc{{\Cal H}_c}
\def\nc{N_{C/X}}
\def\lr{\longrightarrow}
\def\sr{\rightarrow}
\def\lc{L^*\vert_C}
\def\ua{U_{\alpha}}
\def\ub{U_{\beta}}
\def\a{\alpha}
\def\b{\beta}
\def\xa{x_{\alpha}}
\def\xb{x_{\beta}}
\def\ya{y_{\alpha}}
\def\yb{y_{\beta}}
\def\za{z_{\alpha}}
\def\zb{z_{\beta}}
\def\ci{{\Cal I}_{\alpha}}
\def\ca{C_{\alpha}}
\def\p{\partial}
\def\Ya{\left[\dfrac{\p}{\p\ya}\right]}
\def\Za{\left[\dfrac{\p}{\p\za}\right]}
\def\Yb{\left[\dfrac{\p}{\p\yb}\right]}
\def\Zb{\left[\dfrac{\p}{\p\zb}\right]}
\def\c{\gamma}
\def\v{\varphi}
\def\aa{a_{\a}}
\def\ba{b_{\a}}
\def\va{v_{\a}}
\def\noa{N^0_{\a}}

\def\daa{\dfrac{d\aa}{d\xa}}
\def\ta{\theta_{\a}}
\def\tb{\theta_{\b}}
\def\aba{a_{\b\a}}
\def\bba{b_{\b\a}}
\def\cba{c_{\b\a}}
\def\fba{f_{\b\a}}
\def\fab{f_{\a\b}}
\def\fbc{f_{\b\c}}
\def\fca{f_{\c\a}}
\def\aab{a_{\a\b}}
\def\P{\Bbb P}
\def\C{\Bbb C}
\def\mc{{\Cal M}_c}
\def\U{{\Cal C}}
\def\wh{\omega_{\U/\hb}}
\def\sv{s_{\v}}
\def\hv{{\Cal H}_{\v}}
\def\Q{\Bbb Q}

\def\ll{\lambda_1}
\def\lll{\lambda_2}

\document

\heading{\bf $\S$0. Introduction}\endheading

Let $X$ be a $(2n+1)$-dimensional compact complex manifold.
A contact structure on $X$ is
a line subbundle $L$ of $\ct$ such that if $\theta$ is a local section of
$L$, then $\theta\wedge (d\theta)^{n}\neq 0$. This in particular implies
 the canonical bundle $K_X\cong
(n+1) L$. A submanifold $C\subset X$ is called a contact
 submanifold if all  local sections of $L$  vanishes on $C$.  One example of
 complex contact manifold is a odd-dimensional projective spaces.
   The present  work was initiated by
two papers of Bryant. In [Br1] using  Penrose's twistor transform
$\cpc
 \overset T\to\sr S^4$,
he showed
that, among other things,  every
superminimal Riemann surface can be realized by  a smooth  contact curves in
$\cpc$.
 Therefore to study superminimal Riemann surfaces in $S^4$ is the same as
  to study
contact curves in $\cpc$. In another paper, Bryant [Br2] studied relationships
 between manifolds with   ``exotic"
holonomies and  certain contact rational curves on a  complex contact 3-folds.
He showed that
the real slice of the moduli space of  certain contact rational curves
in a complex contact threefold has a so-called ``exotic" ${\Cal
G}_3$-structure.
 Thus
 he found a missing group from  Berger's list of the possible
 pseudo-Riemannian
 holonomies which acts irreducibly on each tangent space.

  The purpose of this paper is to systematically
  study complex contact threefolds and contact curves on them.
   The first main result of this paper  is a classification of
 projective contact threefolds.  Using Mori's theory of
	 extremal ray we   show that the  types of  complex {\it projective}
	 contact threefolds are very limited (Theorem  1.6 and 1.8).
	    The second  result
   is about the moduli space of  contact curves on a complex contact
threefold (Theorem 2.3). We calculate
its Zariski tangent space and the space of obstructions. This result
generalizes a theorem of Bryant [Br2], which dealts with only
rational contact curves.  Finally  we study contact curves in $\cpc$ and
 obtain a Pl\"ucker type formula (Theorem 3.2) for contact curves
 in $\cpc$. It was predicted  in [Br1] that such a formula  should exist.

 The paper is organized as follow.  There will be three sections.
In the beginning of the first section, we will  study some general properties
of a complex contact threefold.  The second half will be devoted to
the classifications of projective contact threefolds. We will start the second
section with an infinitesimal study
of the moduli space of contact curves. We will  relate its normal bundle to
the
 first prolongation of  the dual of  the contact bundle. At the end of the
 second section, we will
 give some applications  of the main theorem in that section.  The last
section  is devoted to studies of  contact curves in $\cpc$.   In that section,
we
have to deal with {\it singular} contact curves.
We  obtain a Pl\"ucker type formula for contact curves in $\cpc$.
Towards the end of that section  we will  relate contact geometry in $\cpc$ to
the geometry of the moduli space of
rank-two stable vector bundles with  $c_1=0$ and $c_2=1$.  We will
 also pose a question on irreducibility of moduli space of
contact curves and speculate how results in the present paper
 can be generalized to higher-dimensional contact manifolds.

  \noindent{\bf Acknowledgments:}  The author wants to express
thanks to his colleague Robert Mclean for very stimulating discussions
 concerning this paper.  He is
 also indebted to his
advisor David Morrison  for his patient teaching and  constant encouragement.

\heading{\bf $\S$1. Compact Complex Contact Threefolds}\endheading

      This section is divided into two parts. In the first part
	  we will study general contact threefolds  even including
non-projective, non-K\"ahler ones. The second part deal specifically with
projective  contact threefolds.

\subheading{\bf General compact complex  contact threefolds}

{}From now on, unless specifically stated, $X$ will be a contact 3-fold
and $L$ is the contact line.  Sometimes
 use a pair $(X,L)$ to mean the same thing.
 In particular  we have $K_X=2L$. Let us
look at some examples of contact threefolds and their contact curves.

\noindent{\bf Example 1.}\hskip 5pt  Let $M$ be a compact complex surface.
 $T_M$ be its
tangent bundle. Then the projectived  bundles
 ${\Bbb P}(T_M)$ a  contact structure.  A word of caution. By ${\Bbb P}(T_M)$
  we mean $T^*_M\setminus \{0\}/{\Bbb C^*}$ instead of $T_M\setminus
\{0\}/{\Bbb C^*}$.  We can write down its local  contact
one forms explicitly. Let $(x,y)$ be a local coordinate system
on $M$ and $\lambda$ be a fiber coordinate. Then $dy+\lambda\,dx$ is
 a local contact one form for ${\Bbb P}(T_M)$. In fact these
local contact forms can be glued together to form a contact line bundle,
which is the dual of the tautological bundle of  ${\Bbb P}(T_M)$.
  The contact curves are either fibers of the
projection of $\pi$, or
horizontal lifting of curves in $M$.  Since $\Omega_M^1\cong
T_M\otimes K_M$, the projectived cotangent bundle
${\Bbb P}(\Omega_M^1)\cong {\Bbb P}(T_M)$. Therefore
${\Bbb P}(\Omega_M^1)$ also have a contact structure. It is
easy to see that its contact line bundle is the dual of the
 tautological line bundle tensored with the pul-back of
 $K_M$.    The following
proposition suggests the contact structure  on ${\Bbb P}(T_M)$
 is  essentially unique.

\proclaim{Proposition 1.1}  For ${\Bbb P}(T_M)$  the moduli space of is
${\Bbb P}\left(GL(T_M)\right)$,
 where $GL(T_M)$ consists of invertible endmorphisms of $T_M$. In particular,
if $T_M$ is simple (for example if $M$ is ${\Bbb C}{\Bbb P}^2$ or a
$K3$ surface), then the contact structure is unique.
\endproclaim
\demo{Proof}  Let
 $X={\Bbb P}(T_M)$ and ${\Cal O}_{X}(1)$ be the tautological line bundle
on $X$. First of all we have the following standard short
exact sequences:
$$0\lr \pi^*\Omega_M^1(1)\lr \Omega^1_X(1)\lr \Omega^1_{X/M}(1)\lr 0\tag
  1.1 $$
$$ 0\lr \Omega^1_{X/M}(1)\lr \pi^*T_M\lr {\Cal O}_X(1) \lr 0 \tag  1.2$$
Since $H^0({\Cal O}_X(1))=H^0(T_M)$, (1.2) implies that $H^0(\Omega^1_{X/M}(1))
=0$.  Sequence  (1.1) in turn gives:
$$H^0(\ct(1))=H^0(\pi^*\Omega_M^1(1))=End\,(T_M)\tag 1.3$$
Now it is clear to see that the space of all contact structures
 is isomorphic to $GL(T_M)$.

\qed
\enddemo

\noindent{\bf Example 2. }\hskip 5pt  The complex projective three space $\cpc$
is
a contact threefold with infinite many contact structures. They are all
equivalent under the actions of the automorphism group of $\cpc$. Their
associated contact line bundles are the same and equal to ${\Cal
O}_{\cpc}(-2)$.  A contact
 structure is a  injective {\it bundle } ({\it not sheaf})  homomorphism
$\varphi$ from ${\Cal O}_{\cpc}(-2)$
to $\pt$. Then it is clear that
$N_{\varphi}=(\text{coker}\,\varphi)^*(-1)$ is a rank-two vector bundle
on $\cpc$ with $c_1=0$ and $c_2=1$. In fact $N_{\varphi}$  is  stable and is a
so-called null-correlation bundle. In [OSS], it is shown
 that any stable rank-two vector bundle on $\cpc$ with
  $(c_1,c_2)=(0,1)$ is a null-correlation
bundle. Moreover two such bundles $N_{\varphi}$ and $N_{\varphi^{\prime}}$
are isomorphic if and only if $\v$ and $\v^{\prime}$ differ by
a non-zero constant. Therefore
the moduli space of contact structures  on $\cpc$ and the space of rank-two
 stable vector
bundles in $\cpc$ with $c_1=0$ and $c_2=1$ are isomorphic. We will come back to
this example at the end of section 3.

 Among all the contact structures, there is a distinguished one obtained
 from the twistor map $\cpc\overset T \to\longrightarrow S^4$. The
contact forms are perpendicular to the tangent directions of the fibers of
$T$.  Bryant [Br1] showed that in the  affine coordinate  chart given by
$[1,z_1,z_2,z_3]$, the local contact form
 is $\omega=dz_1-z_3dz_2+z_2dz_3$.   Curves which are contact with respect to
this contact structure
are called horizontal curves by Bryant. Their images in $S^4$ are so-called
{\it superminimal} surfaces in $S^4$. Therefore to study curves contact with
respect to this contact structure is the same as  to study superminimal
Riemann surfaces in $S^4$.

   By the definition of  contact manifold, to check if a complex 3-fold
$X$ is contact we need to  find a line subbundle $L\subset \ct$
such that  $\theta\wedge d\theta\neq 0$
for any non-zero local section $\theta$ of  $L$. This in particular
implies that $K_X=2L$.  However the  {\it non-integrability} condition
$\theta\wedge
 d\, \theta\neq 0$  is difficult to verify in general. Fortunately, the
following
theorem tells us except for some very special cases  the non-integrability
is satisfied automatically if $K_X=2L$.

\proclaim{Theorem 1.2} Let $X$ be a compact complex threefold and $L\subset
\ct$ such that $K_X=L^2$. If  $L$ does not define a contact structure
on $X$, then  there is a  smooth fibration $X\overset \v\to \lr C$ such that
 \roster
\item $C$ is a smooth  curve  and $L=\v^*(K_C)$.

\item $\v$ has connected fibers and they  all have trivial canonical
 bundles.
\endroster
If moreover $X$ is K\"ahler, then $X$ is not simply connected and
fibers of $\v$ are either
$K3$ surfaces or abelian surfaces.
\endproclaim

\demo{Proof}   Choose a local coordinate cover $\{\ua\}_{\a\in\Gamma}$ on $X$
such that $L\vert_{\ua}={\Cal O}_{\ua}(\ta)$, where $\ta$ is a local
holomorphic one form. Consider the ${\Cal O}_{\ua}$-module
homomorphism $\psi: L^2\sr K_X$ defined by $\psi(\ta^{\otimes 2})=\ta\wedge
d\ta$.
It is clear that $\psi$ is well-defined. Since $K_X=L^{2}$, $\psi$ is either
 nowhere zero or identically zero. If $\psi$ is nowhere zero, then
$L$ defines a contact structure. If $\psi$ is identically zero, then
$\ta\wedge d\ta=0$ for any $\a\in \Gamma$. We are going to show that in this
 case  $X$ admits a fibration described above.

 For any $\a\in \Gamma$ let $(\xa,\ya,\za)$ be a coordinate system. Set
$\ta=a_{\a}d\xa+b_{\a}d\ya+c_{\a}d\za$. Then the fact that
 $\ta\wedge d\ta=0$ implies that $a_{\a},b_{\a}$  and $c_{\a}$ are constant
 proportional to one another. Hence $\ta=g_{\a}(\lambda_{\a}d\xa+
\mu_{\a}d\ya+\sigma_{\a}d\za)$, where $g_{\a}$ is a holomorphic function
on $\ua$ and $\lambda_{\a},\mu_{\a}$ and $\sigma_{\a}$ are constant.
Since $\ta$ generates $L$ over $\ua$, the function
$g_{\a}$ has to be invertible in ${\Cal O}_{\ua}$. Therefore, we may well
assume $g_{\a}=1$. Then we can write $\ta=df_{\a}$, where
$f_{\a}=\lambda_{\a}\xa+
\mu_{\a}\ya+\sigma_{\a}\za$.  Hence we can choose coordinates in such a
 way such that $\ta=d\xa$ for each $\a\in\Gamma$.

Let $\aab$ be the transition function for $L$ over $\ua\cap\ub$. Then
On $\ua\cap\ub$, we have $\ta=\aab\tb$, i.e., $d\xa=\aab d\xb$.  This implies
that $\dfrac{\p\xa}{\p\yb}=\dfrac{\p\xa}{\p\zb}=0$, i.e.,
 $\xa$ depends on $\xb$ only. Let $\xa=\fab(\xb)$ for some
one-variable holomorphic function $\fab$. It is clear that
$\fab$ satisfies the co-cycle conditions, i.e., $\fab\circ\fba=id.$ on
$\ua\cap\ub$ and $\fab\circ\fbc\circ\fca=id.$ on $\ua\cap\ub\cap U_{\gamma}$.
Then $\{\fab\}$'s define a complex curve $C$. It is clear that $C$ is smooth
and compact. There
is a natural holomorphic map $X\overset \v\to \lr C$ sending $(\xa,\ya,\za)$ to
 $\xa$. Since $\aab=\dfrac{d\fab}{d\xb}$, we have $L=\v^*(K_C)$. The inclusion
 $L\subset \ct$ fits into the natural short exact sequence:
$$0\lr L=\v^*K_C\overset d\v^*\to \lr \ct \lr \Omega^1_{X/C}\lr 0\tag 1.4$$
 where $\Omega^1_{X/C}$ is the relative cotangent sheaf. The fact that
$d\v^*$ is a bundle injection implies that $\Omega_{X/C}^1$ is locally
free, hence $\v$ is a smooth fibration. By passing to the Stein
factorization, we can assume that $\v$ has connected fibers and
$\v_*({\Cal O}_X)= {\Cal O}_C$.  Let $K_{X/C}=\wedge^2\Omega_{X/C}^1$ be the
 relative canonical
bundle. Then $K_{X/C}=\v^*(K_C)$ since $K_X=2L=2\v^*(K_C)$. Therefore
$K_F$ is trivial for any fiber $F$ of $\v$.

If moreover $X$ is K\"ahler, then any fiber $F$ is also K\"ahler.
Therefore $F$ is
either a $K3$ surface or an Abelian surface.   Since $K_{X/C}=\v^*(K_C)$,
 we have $\v_*K_{X/C}=K_C$. By a theorem of Fujita-Kawamata [Ka1]
$\v_*K_{X/C}$ is weakly
positive. Hence $C$ is irrational.  Leray spectral sequence
 implies that $h^1({\Cal O}_X)=h^0(R^1\v_*{\Cal O}_X)+h^1({\Cal O}_C)$.
Since $C$ is irrational, $h^1({\Cal O}_C)>0$. Hence $h^1({\Cal O}_X)
>0$, therefore $b_1(X)=2h^1({\Cal O}_X)>0$. Hence $X$ is not simply connected.
\qed\enddemo

  The above theorem implies the following corollary immediately.
\proclaim{Corollary 1.3} If $X$ is a simply connected K\"ahler threefold and
 $L\subset \ct$ such that $K_X=2L$, then $L$ defines a contact structure.
\endproclaim

\proclaim{Theorem 1.4} Let $X$ be a compact K\"ahler threefold with $c_1(X)=0$.
Then $X$ has no contact structure.
\endproclaim
\demo{Proof}  By Beauville's theorem  [Be],
there is a finite unramified cover $\tilde X$ of $X$ such that
$\tilde X=T\times Y$, where $T$ is  a complex torus  and $Y$ is
   simply connected  with trivial canonical bundle.
This in particular implies that $K_X$ is a torsion line bundle.
Suppose $X$ has a contact structure with a contact line bundle
$L$. We distinguish two cases:

\noindent Case I: $\text{dim}\,Y>0$.
 Since $K_X$ is torsion and $K_X=L^2$, we
see that $L$ is also torsion. Let $m$ be a non-negative integer such that $L^m$
is
trivial. On the one hand,  the bundle injection $L\hookrightarrow \ct$ implies
a bundle injection ${\Cal O}_X=L^m\hookrightarrow \ct$. Hence $h^0(
S^m\ct)>0$. Since $\tilde X$ is an unramified cover,
$h^0(S^m\Omega_{\tilde X}^1)>0$.  On the other hand,
$h^0(S^m\Omega_{\tilde X}^1)=\sum_{p+q=m}h^0(S^p\Omega_T^1)\,
h^0(S^q\Omega_Y^1)$. A theorem of Kobayashi [Ko] implies
that $h^0(S^q\Omega_Y^1)=0$ for all $q\geq 0$. Hence
 $h^0(S^m\Omega_{\tilde X}^1)=0$. This gives a contradiction.

\noindent Case II:  $\text{dim}\,Y=0$. In this case, $\tilde X$ is a
complex torus. The  given contact structure $L$ on $X$ induces a contact
structure $\tilde L$ on $\tilde X$. Since $\Omega_{\tilde X}^1$ is trivial,
 $h^0(\tilde L^{-1})>0$. However $c_1(\tilde L)=0$. Therefore
$\tilde L$ has to be trivial. But a trivial line bundle
on a torus does not define a contact bundle. This can be shown as follow. Note
 that  $\tilde L$ gives rise to
a global  holomorphic one form $\theta$ on $\tilde X$. Let
$f: {\Bbb C}^3\lr \tilde X$ be the universal cover, and $\omega=f^*\tilde
\theta$. Then $\omega$ has to be of the form
$adx+bdy+cdz$, where $a,b,c$ descend to holomorphic
functions on $\tilde X$. Hence they must all be constant. Hence $d\omega=0$.
Therefore $\tilde\theta\wedge d\tilde\theta=0$. This means that
$\tilde \theta$ does not define a contact structure.

  In either case we get a contradiction. Hence we are done.

\qed
\enddemo

Next we will give a topological obstruction
  to existence of a contact structure.

\proclaim{Proposition 1.5} Let $X$ be a compact complex contact threefold with
a contact
line bundle $L\subset \ct$. Then
 $$\chi_{top.}(X)=12\chi({\Cal O}_X)-\dfrac{c_1(X)^3}{8}$$
  where $c_i(X)=c_i(T_X)\, \,
 i=1,2,3$, and
$\chi({\Cal O}_X)$ is the holomorphic Euler characteristic, and
$\chi_{top.}(X)$ be
 the topological Euler characteristic.

\endproclaim

\demo{Proof}  Since $L\subset \ct$ is
a subbundle (rather that a subsheaf),  we conclude that \linebreak
$c_3(\ct\otimes L^{-1})=0$
by Porteous' formula. This implies that:
$$-c_3-c_2L-c_1L^2-L^3=0 \tag 1.5 $$
where $c_i=c_i(X)$, for $i=1,2,3$. Using the fact that $c_1=-K_X=-2L$ and
Riemann-Roch formula $\chi({\Cal O}_X)=\dfrac{c_2c_1}{24}$, we get:
$$c_3=12\chi({\Cal O}_X)-\dfrac{c_1(X)^3}{8}\tag 1.6$$
 The fact that $c_3=\chi_{top.}(X)$ implies the proposition immediately.

\qed
\enddemo

\subheading{\bf   Projective contact threefolds}

  In this subsection, we will assume that $(X,L)$ is a  complex projective
contact
threefold. By this we mean that $X$ is a projective
 complex manifold and $L\subset \ct$ defines a contact structure
on $X$.  We will show that  the types of these contact threefolds are very
limited.

 Before we state our next theorem, let us recall that a line bundle
 is called {\it nef} if its intersection with every  effective curve is
non-negative.

\proclaim{Theorem 1.6} If $(X,L)$ is a projective complex  contact threefold
and $K_X$ is {\it not} nef, then $X$ is either isomorphic to  $\cpc$ or
$(X,L)\cong \left({\Bbb P}(T_M),
{\Cal O}_{{\Bbb P}(T_M)}(-1)\right)$ for some
smooth complex projective surface $M$. \endproclaim
\demo{Proof} Since $K_X$ is not nef, Mori's theory of
extremal ray implies that there is a smooth rational curve $C\subset X$
such that $4\geq -K_X\cdot C>0$ and $C$ generates an extremal ray
$R={\Bbb R}_{+}[C]$. If $-K_X\cdot C=4$, then it is well-known that  $X$ is
isomorphic to $\cpc$. Otherwise $-K_X\cdot C=2$ since $K_X=2L$.
This implies that $L\cdot C=-1$.    Consider the restriction
of the contact sequence to $C$:

$$
\CD
@. @. 0 @.\\
@. @.  @AAA  @.\\
@. @.\nc  @.\\
@.  @. @AAA  @.\\
0 @>>> L^{\bot}\vert_C @>>> T_X\vert_C @>>> L^*\vert_C @>>> 0 \\
@. @.  @AAA @.\\
@. @.  T_C  @.\\
@. @. @AAA @.\\
@. @.  0 @.
\endCD \tag 1.7
$$
 Denote by  $\alpha$  the composed map
$ T_C\sr \lc$. Since $T_C\cong {\Cal O}_C(2)$ and
$\lc\cong {\Cal O}_C(1)$, we conclude that  $\alpha$ has to be zero, i.e., $C$
is
 a contact curve.
Then by a theorem of Bryant [Br2] or Theorem 2.3 below,
$\nc\cong {\Cal O}_C\oplus {\Cal O}_C$. It is shown by Mori [Mo] that
$X$ is isomorphic to a conic bundle in this case.  Since $L\cdot C=-1$, any
deformation
of $C$ in $X$ is still reduced and irreducible. Hence $X$
 must be a ${\Bbb C}{\Bbb P}^1$-bundle over some smooth projective
 surface $M$.  Therefore we can write $X$ as ${\Bbb P}(E)$ for some
rank-two vector bundle on $M$. Let $X\overset \pi \to\lr M$ be the natural
projection. Note that any fiber of $\pi$ generates the extremal ray $R$.  Let
${\Cal O}_{X}(1)$ be the tautological line bundle of $X$.
It is easy to see that $K_X=-2{\Cal O}_{X}(1)+
\pi^*\left(K_M+\wedge^2 E\right)$. However $K_X=2L$.
Hence $\pi^*\left(K_M+\wedge^2 E\right)=-2L_0$, where $L_0=
{\Cal O}_{X}(1)-L$.  Since $L_0\cdot C=-1+1=0$,
$L_0=\pi^*L_1$ for some line bundle $L_1$ on $M$. Therefore
$K_M+\wedge^2 E=-2L_1$, i.e., $\wedge^2\left(E\otimes L_1\right)=\wedge^2T_M$.
Hence
if we tensor $E$ by   $L_1$, then
we can  assume that $\wedge^2 E\cong \wedge^2 T_M$
and  ${\Cal O}_{X}(-1)$ is  the contact line bundle.
 Then
there is a natural  bundle injection $\lambda: {\Cal O}_{X}(-1)\lr
 \ct$.

We will show that $E\cong T_M$. We first prove the following claim:

\noindent{\it Claim:} Let $\cmt$ be the relative cotangent bundle. Then
$$H^0\left(\cmt(1)\right)=0,\,\,\, H^1\left(\cmt\right)\cong \Bbb C $$

\demo{Proof of the Claim}  Consider the relative Euler sequence:
$$0\lr \cmt(1)\lr \pi^*E\lr {\Cal O}_{X}(1)\lr 0 \tag 1.8$$
Since $H^0(\pi^*E)\cong H^0(E)\cong
H^0\left({\Cal O}_{X}(1)\right)$, we have
$H^0\left(\cmt(1)\right)=0$.
Since $\pi$ is a ${\Bbb C}{\Bbb P}^1$-bundle,
we have $\pi_*\cmt=0$. Hence Leray spectral sequence for $\pi$
implies that $\h1\cmt\cong H^0\left(R^1\pi_{*}\cmt\right)$.
However  by the  relative duality, $R^1\pi_{*}\cmt\cong \left(\pi_*{\Cal O}_X
\right)^*
\cong {\Cal O}_M$. Hence
$\h1\cmt\cong \Bbb C$. Hence the claim is proved.

\enddemo

   Now consider the tangential sequence:
$$ 0\lr T_{X/M}\lr T_X\lr \pi^*T_M\lr 0 \tag 1.9$$
Since $H^0\left(\cmt(1)\right)=0$, the  bundle
injection (from the contact structure on $X$) $\lambda:{\Cal O}_{X}(-1)\lr
 \ct$ induces a surjective bundle map
$\sigma: \pi^*T_M\lr {\Cal O}_{X}(1)$. Let ${\Cal N}$ be the
kernel of $\sigma$. Then ${\Cal N}$ is a line bundle. Moreover
$\Cal N\cong \pi^*\left(\wedge^2 T_M\right)\otimes
{\Cal O}_{X}(-1)
$. However by our assumption, $\wedge^2 E\cong \wedge^2 T_M$. Therefore
 sequence (1.8)  implies that ${\Cal N}\cong \cmt(1)$. Hence
 we obtained an exact sequence:
$$0\lr \cmt(1)\lr \pi^*T_M\lr {\Cal O}_{X}(1)\lr 0 \tag 1.10$$
 Let $e_1$, respectively $e_2$ be  the extension class corresponding to
 (1.8), respectively (1.10). They are elements in $\h1\cmt\cong \Bbb C$.
They are non-zero  since their
corresponding sequences do not split (because their restrictions
to a fiber of $\pi$  do not split since the normal
bundle of the fiber is trivial).
Hence  they differ only by a non-zero scalar.
 This last fact implies that $\pi^*E
\cong \pi^*T_M$. Since $\pi_*{\Cal O}_X={\Cal O}_M$, we have $E\cong T_M$.
 Hence we are done.

\qed
\enddemo

\proclaim{Corollary 1.7} If  $X$ is Fano, i.e., $-K_X$ is
ample, then $X$ is isomorphic  either to $\cpc$ or
$\P\left(\Omega^1_{\cpb}\right)$.
\endproclaim
\demo{Proof} Suppose that $\cpc$ or $(X,L)\cong \left({\Bbb P}(T_M),
{\Cal O}_{{\Bbb P}(T_M)}(-1)\right)$ for some
smooth complex projective surface $M$. Then
$-K_X={\Cal O}_{{\Bbb P}(T_M)}(-1)$. Hence by our assumption, $T_M$
is ample. This clearly implies that $M\cong \cpb$ by Mori's proof
of Hartshorne conjecture. However as see in Example 1, ${\Bbb P}(T_{\cpb})\cong
\ppb$.
Hence we are done.

\qed
\enddemo

\noindent {\bf Remarks:} We say that  two contact manifolds $X_1$ and
$X_2$ are {\it contactly birational} to each other if there is a
 birational map between them such that contact curves are mapped to
 contact curves. It is showed in [Br1] that $\cpc$ and $\ppb$ are
 contactly birational to each other. The above Corollary 1.7 answers
affirmatively  a
 question of Bryant, who asked if any contact threefold with positive first
 Chern class is contactly birational   to $\cpc$.

  Now we are going to study projective contact threefolds with
nef canonical bundle. They are  so-called minimal model. We will show that
they are either fibered by abelian surfaces, or hyperelliptic surfaces
 or elliptic curves.

\proclaim{Theorem 1.8} Let $(X,L)$ be a projective  complex contact threefold
such
that $K_X$ is nef.
Then the Kodaira dimension of $X$ is either one or two. If its Kodaira
 dimension is one, then $X$ admits a fibration  onto a smooth
curve such that its generic fiber is an abelian surface or, a
hyperelliptic surface. If its
Kodaira dimension is two, then it admits an elliptic fibrations.
\endproclaim

\demo{Proof} Let $\kappa(X)$ be the Kodaira dimension of $X$. Since
$K_X$ is nef, then $\kappa(X)$ is between 0 and 3. If $\kappa(X)=0$, then
$c_1(X)=0$. We have shown  in Theorem 1.4 that  this can not
happen. Hence $\kappa(X)=1,2$, or 3.  We are going to
distinguish three cases according to $\kappa(X)$.

\noindent Case 1: $\kappa(X)=3$. In this case $K_X$ is big, i.e., $K_X^3>0$.
Then $L^3>0$ also. By a theorem of Tsuji [Ts], $T_X$ (hence $\ct$ too) is
$K_X$-semistable. Since $L\subset \ct$ is a subbundle, we have
$L\cdot K_X^2< \dfrac{K_X\cdot K_X^2}{3}$, i.e., $4L^3<\dfrac{8}{3}L^3$. This
 is absurd. Hence $\kappa(X)$ can not be 3.

\noindent Case 2: $\kappa(X)=2$. Then by abundance theorem [Ka2], for a
sufficiently large $m$, $|mK_X|$ is free.  Let $f$ be the morphism defined
 by $|mK_X|$ for some fixed large integer $m$. Let $S$ be the image of $f$.
Then $S$ is a normal complex
surface. By passing to its Stein factorization, we can assume that
$f:X\lr S$ has connected fibers and $f_*{\Cal O}_X={\Cal O}_S$. By the
definition
 of $f$, we see that there is an ample line bundle $H$ over $S$ such that
$mK_X=f^*H$. Let $F$ be a generic fiber of $f$. Then $F$ is smooth curve and
 by adjunction formula $K_F=K_X\vert_F$. Since
 $mK_X\vert_F\cong {\Cal O}_F$ Therefore $K_F$ is torsion. This implies that
$F$ is
an elliptic curve.

\noindent Case 3: $\kappa(X)=1$. By abundance theorem  [Mi] again, we know
 that  $|mK_X|$  is free for sufficiently large $m$.  Let $f:X\lr C$ be the
morphism defined by $|mK_X|$. Then $C$ is a normal
 projective curve. Hence $C$ is smooth. Let $F$ a general fiber of
$f$. Then $F$ is a smooth projective surface.  As we did in the previous
case, we  can show that $mK_F\cong {\Cal O}_F$, i.e., $K_F$ is
a torsion line bundle. The rest part of the proof
is very much similar to the proof of
Theorem 1.4 above.  Then $F$ either has finite foundmental
group, or its universal cover is a complex torus.  We claim that
$F$ can not have finite foundmental group. This can be shown as follows.
The bundle injection $L\hookrightarrow\ct$ induces  a
bundle injection $L\vert_F\hookrightarrow \Omega_F^1$. Since
$2mL\vert_F\cong mK_X\vert_F\cong {\Cal O}_F$, we get
a bundle injection ${\Cal O}_F\hookrightarrow S^{2m}\Omega_F^1$. Hence
$H^0\left(S^{2m}\ct\right)>0$. This contradicts  a theorem of Kobayashi [Ko].
Therefore $F$ must have infinite foundmental group. Hence $F$ is
either  an abelian or  a hyperelliptic surface.

\qed

\enddemo

  The above theorem gives rise to a natural question:  do those
manifolds described in the theorem do have a contact structure?. My guess
 is that to be a contact manifold the fibration should at least have no
singular fibers.
   Next we show that the trivial abelian fibration  over
   $\C\P^1$ does have a
contact structure.

\noindent{\bf Example 3.} \hskip 5pt  Let $A$ be any
abelian surface.  Let $X=\C\P^1\times A$. Let $(x,y)$ be the
coordinates on $\Bbb C^2$, which is the universal cover of $A$. Let
$U_1$ and $U_2$ be the two coordinate cover of $\C\P^1$, and $z_1$ and
$z_2$ be respective coordinates. Then over $U_1\cap U_2$, $z_1=1/z_2$.
Consider two local   one-forms $\theta_1=dy+z_1\, dx$ on $U_1\times A$, and
$\theta_2=z_2\, dy+dx$ on $U_2\times A$.  Then these two local forms define a
contact structure on
$X$. In fact the contact line bundle $L$ is just $p_1^*{\Cal O}_{\C\P^1}(-1)$,
where $p_1:X\sr \C\P^1$ is the first projection.

\heading{\bf $\S$2. Contact Curves and Their Moduli Space }\endheading

 Let $X$ be a  compact complex contact threefold with a contact line bundle
$L\subset \ct$. Let $C$ be a smooth contact curve in $X$. $\hb$ be
the irreducible component of the Douady space of $X$ the contains $[C]$.
Let $\hc\subset \hb$ be the set of all contact curves. The purpose of this
section is to study the subspace $\hc$.

\proclaim{Lemma 2.1}[Bryant] Let $N_{C/X}$ be the normal bundle of $C$. Then
 we have the following short exact sequence:
$$0\lr L^*\otimes K_C \lr \nc\lr \lc\lr 0\tag 2.1$$
\endproclaim
\demo{Proof}
 Consider the  sequence (1.7).
Let $\alpha: T_C\lr L^*\vert_C$ be the composed homomorphism. Since $C$ is
a contact curve, $\alpha=0$. This induces a surjective bundle map
 $\nc\overset \beta\to \lr L^*\vert_C\lr 0$. Let $F$ be the kernel of $\beta$.
 Then $F$ is a line bundle. By the  adjunction formula, we have:
$$\wedge^2\nc=K_C\otimes K_X^*=K_C\otimes L^{*2} \tag 2.2$$
However it is clear that $\wedge^2\nc=F\otimes L^*\vert_C$. Now (2.2)
implies that $F=K_C\otimes L^*$. Therefore the lemma is proved.
\qed
\enddemo

Before we continue, let us  recall the definition and some properties
of the {\it first prolongation}  of a line bundle $N$ on a compact
complex manifold, say $Y$. The first prolongation of $N$, denoted by
$P^1(N)$,  is a rank-two bundle obtained via the following exact sequence:
$$0\lr \Omega_Y^1(N)\lr P^1(N) \lr N\lr 0\tag 2.3$$
where the extension class of (2.3) is   $c_1(N)\in
H^1(\Omega_Y^1)=\text{Ext}_Y^1(N,\Omega_Y^1(N))$. If $\{g_{\a\b}\}$ are
 transition functions for $N$, then the Cech co-cycle
$\{d\log(g_{\a\b})\}$ represents $c_1(N)$, hence the extension class for (2.3).

\proclaim{Lemma 2.2}
\roster
\item Sequence (2.3) splits  if and  only if  $c_1(N)=0$, i.e., $N$ is flat.
\item  Sequence (2.3) always splits cohomologically, i.e.,
$H^0(P^1(N))\sr H^0(N)$ is always surjective.
\endroster
  \endproclaim

\demo{Proof} The first part is true by definition. Let us prove the second
 part.  Choose a local coordinate cover $\{\ua\}$ on $Y$ such that
$N$ is trivialized on each $\ua$. Let $\{g_{\a\b}\}$'s be the transition
 functions
for $N$ under these trivializations. Let $s=\{s_{\a}\}$ be an
 arbitrary global holomorphic section of $N$. The second part is equivalent to
show that the
 Cech co-cycle $\{s_{\a}d\log(g_{\a\b})\}$ is a co-boundary with
coefficient in $\Omega_Y^1(N)$, i.e., it is zero in $H^1(\Omega_Y^1(N))$.
 Since $s$ is
a global section of $N$, $s_{\a}=g_{\a\b}s_{\b}$. This implies that
 $s_{\a}d\log(g_{\a\b})=ds_{\a}-g_{\a\b}ds_{\b}$. Hence
$\{s_{\a}d\log(g_{\a\b})\}$ is a co-boundary. We are done.

\qed
\enddemo

We now state our next theorem.

\proclaim{Theorem 2.3} Let $C\subset X$ be a smooth contact curve in
a contact threefold. Let $L\subset \ct$ be the contact structure. Then
 the following are true:
\roster
\item  The Zariski tangent space $T_{[C]}\hc\cong H^0(L^*\vert_C)$. The space
 of obstructions is $H^1(L^*\vert_C)$. In particular, if $h^1(\lc)=0$, then
$\hc$ is smooth of dimension $h^0(\lc)$.
\item  The normal bundle $\nc$ is isomorphic to $P^1(\lc)$, where $P^1(\lc)$ is
the
first prolongation  of $\lc$. In particular,  sequence (2.1)
splits if and only if $L\cdot C=0$.
\endroster
\endproclaim

\demo{Proof} As Bryant showed in [Br2] that we can choose a local coordinate
 cover $\{U_{\alpha}\}$on $X$ such that on each $\ua$, $L$ is generated
by a local one form $\theta_{\alpha}$ of the form
$\theta_{\alpha}=dy_{\a}-z_{\a}dx_{\a}$, where $(\xa,\ya,\za)$ are local
coordinates and
$C\cap\ua=\{\ya=\za=0\}$. Denote $C\cap\ua$ by $C_{\a}$. Let $\Ya$ and  $\Za$
be the classes of
$\dfrac{\p}{\p\ya}$ and $\dfrac{\p}{\p\za}$ in $\nc\vert_{\ca}$. Then they
generate $\nc\vert_{\ca}$ over ${\Cal O}_{\ca}$.

 We call an embedded deformation $\{C_t\}$ of  $C$ is a
 {\it contact} deformation if
$\forall t$, $C_t$ is a contact curve with respect to
the given contact structure.  Given a  normal vector field
 $v_{\a}=a_{\a}\Ya+b_{\a}\Za$ with $a_{\a}$ and $b_{\a}$ being in ${\Cal
O}_{\ua}$.
Then $v_{\a}$ generates a first-order contact deformation of $\ca$ if and only
if
$b_{\a}=\dfrac{da_{\a}}{d\xa}$. This is because if  $\va$ generates a
one-dimensional deformation  $C_{\a t}$ of $\ca$ in $\ua$, then  the
deformation  can be expressed
as:
$\xa(t)=\xa+c_{\a}t+O(t^2), \ya(t)=\aa t+O(t^2), \za(t)=\ba t+O(t^2)$.
Then $\theta_{\a}\vert_{C_{\a t}}=t\,(\dfrac{da_{\a}}{d\xa}-\ba)d\xa +
O(t^2)$.
Therefore $\va$ generates an infinitesimal contact deformation of $\ca$
if and only if $b_{\a}=\dfrac{da_{\a}}{d\xa}$.

  Let us denote by $\noa=\bigg\{\aa\Ya+ \daa\Za\,\bigg\vert\, \aa\in {\Cal
O}_{\ca}
\bigg\}$.
Then $\noa$ is a subset of $\nc\vert_{\ca}$, but not a submodule over
${\Cal O}_{\ca}$. However, we can put a different   ${\Cal O}_{\ca}$-module
structure on $\noa$.
For any $f_{\a}\in {\Cal O}_{\ca}$, we define its action on
$\aa\Ya+ \daa\Za$ by:
$$f_{\a}\circ \left(\aa\Ya+ \daa\Za\right)\overset\text{def.}\to =
\left(f_{\a}\aa\right)\Ya+
\dfrac{d\left(f_{\a}\aa\right)}{d\xa}\Za $$
It is clear that with this ${\Cal O}_{\ca}$-module structure $\noa$ is free and
generated by $\Ya$. Next we will show that we can glue these $\noa$'s to
get a line bundle over $C$. In fact  we will
show that this line bundle is isomorphic to
$\lc$. First let us prove the following claim.

\noindent{\it Claim :  Let $\ci=(\ya,\za)$ be the ideal sheaf of $\ca$.
 After modulating  the ideal sheaf $\ci^2$,  on $\ua\cap\ua$ we have:}
$$
\left\{
\aligned
\xb &= \fba(\xa)\\
\yb &= \aba\ya\\
\zb &= \cba\za+\bba\ya
\endaligned \right.
$$
{\it where $\fba$, $\aba$, $\bba$ and $\cba$ are  holomorphic
functions over $\ua$ depending only on $\xa$.
Moreover $\bba=\dfrac{d\aba}{d\xb}$ and
$\cba=\aba\dfrac{d\xa}{d\xb}$ when restricted to $\ua\cap\ub\cap C$. }

\demo{Proof of  the Claim }  On $\ua\cap\ua$, we have:
$$d\xb=\dfrac{\p\xb}{\p\xa}d\xa+\dfrac{\p\xb}{\p\ya}d\ya+
\dfrac{\p\xb}{\p\za}d\za  \tag 2.4$$
Note that $d\xa$ generates $\Omega_C^1\vert_{\ca}$, therefore
$d\xb$ and $d\xa$ are proportional to each other on $\ua\cap\ub\cap C$.
Therefore
$$\dfrac{\p\xb}{\p\ya}=\dfrac{\p\xb}{\p\za}=0 \,\,\,\,(\text{mod}\, \ci)
\tag 2.5
$$
 This implies the first equation of the claim.

By the same token, since
 $\ta$ and $\tb$ are proportional, we have
$\dfrac{\p\yb}{\p\za}-\zb\dfrac{\p\xb}{\p\za}=0$. Therefore (2.5) implies that
$\dfrac{\p\yb}{\p\za}=0\, (\text{mod} \, \ci^2)$. Hence
$\yb=h_{\a\b}(\ya,\xa)\,(\text{mod}\, \ci^2)$. Since $\ya$ and $\yb$ are
both zero when restricted to $\ua\cap\ub\cap C$, we see $\yb$ has to
 be divisible by $\ya$. This give the second equation of the claim. The last
equation of  the claim is obvious since $\zb$ vanishes on $\ua\cap\ub\cap C$.

The rest of the claim follows easily  from the fact that $\ta$ and $\tb$ are
proportional to each other on $\ua\cap\ub$. Hence the claim is proved.

\enddemo

On the one hand, by the claim, when restricted to $\ua\cap\ub\cap$ C, we have:
$$
\aligned
\Ya &= \dfrac{\p\yb}{\p\ya}\Yb+\dfrac{\p\zb}{\p\ya}\Zb\\
   &=  \aba\Yb+\dfrac{d\aba}{d\xb}\Zb\\
 &= \aba\circ \Yb
\endaligned \tag 2.6$$

On the other hand, it is clear to see that $\tb=\aba\ta\, (\text{mod}\, \ci)$
 on $\ua\cap\ub\cap C$. Therefore we conclude that $\noa$'s can be glued to
get a line bundle on $C$ which is  isomorphic to $\lc$. This proves the
 first part of the
theorem.

Now let us prove the second part of the theorem.  By the claim, the
transition matrix for the normal bundle $\nc$ over $\ua\cap\ub\cap C$ is
given by:
$$
\aligned
\Ya &= \aba\Yb+\dfrac{d\aba}{d\xb}\Zb\\
\Za &= \aba \dfrac{d\xa}{d\xb} \Zb
\endaligned\tag 2.7
$$

 This implies that the extension class (2.1) is given by $e_{\a\b}=
\aab\dfrac{d\aba}{d\xb}$. Therefore $e_{\a\b}d\xb=\aab\,d\aba=d\log
(\aba)$ gives a class in $H^1(K_C)$, which is the extension class
 corresponds to the exact sequence (2.1). Hence the normal
bundle $\nc\cong P^1(\lc)$ by Lemma 2.2.  Lemmas 2.2 also implies that
sequence (2.1) splits if and only  if \linebreak $c_1(\lc)=0$, i.e., $L\cdot
C=0$.
\qed
\enddemo

Next we will study the moduli space of all curves which are  contact with
respect to
 a given contact structure on $X$.  Let $\mc$ be the moduli space of all
contact structures on $X$.  Let $\v\in\mc$ be a contact structure on $X$ and
 $L_{\v}$ be the contact line bundle. Since
$K_X=2L_{\v}$, if Pic($X$) has no two torsion, then
all the $L_{\v}$ are the same.  To give   an element $\v\in \mc$ is the
 same as to give  a bundle injection ( which is still denoted by $\v$)
 $L_{\v}\overset \v\to\lr \ct$. Therefore $\mc$ is an Zariski open subset
of $H^0\left(\Omega_X^1(L_{\v}^{*})\right)$.

  Let $\hb$ be an irreducible component
of the Douady space of {\it smooth} curves in $X$. Let $\U\subset X\times \hb$
be the universal
family. The we have two projections, namely,  $\U\overset p\to\lr \hb$ and
 $\U\overset q \to \lr X$. Let $\wh$ be the relative dualizing sheave for
$p$.  Given a contact structure $\v$  on $X$, it induces in a  natural  way a
section
$s_{\v}$ of  $\E_{\v}\overset \text{def.}\to=p_*\left(\wh\otimes
q^*L_{\v}^*\right)$. The
vanishing locus of $s_{\v}$ is $\hv$, the set of curves in $\hb$ contact
with respect to $\v$. This proves the first part of the following proposition.

\proclaim{Proposition 2.4} Let $X$ be a compact complex threefold. Then the
following is true:
\roster
\item Given a contact structure $\v$ on $X$ there is a natural section
$\sv$ of $\E_{\v}$ such that $\hv$ is  the vanishing locus of $s_{\v}$.
\item   If $H^0(L_{\v}\vert_C)=0$ and  $\hv$ is non-empty, then the smoothness
of
$\hv$ at a point $[C]\in \hv$ is equivalent to the smoothness of $\hb$ at
the same point.
\item Let $[C]\in \hv$ be contact curve. If
$H^1\left(L_{\v}^*\vert_C\right)=0$, then
both $\hb$ and $\hv$ is smooth in a neighborood  of $[C]$.
\item If $\hv$ is smooth everywhere, then
$ T_{\hb}\vert_{\hv}\cong T_{\hv}\oplus \E_{\v}$, i.e., the restriction of
tangent bundle splits. The normal bundle  of $\hv$ in $\hb$ is isomorphic to
 $\E_{\v}$.
\endroster
\endproclaim
\demo{Proof} The first part is proved in the above paragraph. As for
the second part,  Lemma 2.2 and Theorem 2.3 imply that  sequence (2.1)
splits cohomologically. Hence the following sequence is exact:
$$
0\lr H^1\left(L_{\v}^*\otimes K_C\right) \lr H^1\left(\nc\right)
\lr H^1\left(L^*_{\v}\right) \lr 0\tag 2.8 $$
Since by our assumptions $H^0\left(L\vert_C\right)=0$, we have
$H^1\left(L_{\v}^*\otimes K_C\right)=0$ by Riemann-Roch theorem. Hence
$H^1\left(\nc\right)
\cong H^1\left(L^*_{\v}\right)$. By local computations as we did in the proof
of Theorem 2.3, we can show that obstructions in $H^1\left(\nc\right)$ are
mapped to obstructions in $H^1\left(L^*_{\v}\right)$.
 Therefore if  $\hv$ is smooth at $[C]$, then
  $\hb$ is also  smooth at $[C]$. Conversely, if $\hb$ is smooth
   at a point $[C]$, then its dimension around that point is
   $h^0\left(\nc\right)$. The fact that $\hv$ is the
   vanishing locus of  the section $s_{\v}\in H^0\left(\E_{\v}\right)$ implies
that
 the dimension of $\hv$  around the given point
$[C]$ is at least $h^0\left(\nc\right)-h^0\left(L^*_{\v}\otimes K_C\right)=
h^0\left(L^*_{\v}\right)$. By Theorem 2.3,
$T_{[C]}\hv=H^0\left(L^*_{\v}\right)$. Therefore $\hv$ is smooth
 at $[C]$. Hence the smoothness
of
$\hv$ at a point $[C]\in \hv$ is equivalent to the smoothness of $\hb$ at
the same point.  This proves the part two.

In view of the fact that sequence (2.1) splits cohomologically, the last two
 parts are clear.

\qed

\enddemo

Theorem 2.3 and Proposition 2.4 imply that  the following corollary,  which was
proved by Bryant [Br2] using
 a different method.
 \proclaim{Corollary 2.5} [Bryant] \hskip 3pt  Let $(X,L)$ be a contact
threefold. Suppose
 $C\subset X$ is a   smooth contact  rational curve in $X$ with
 $L\vert_C\cong {\Cal O}_C(-k-1)$ for some integer $k\geq 0$. Then the moduli
space of
 contact curves that contains $C$ is smooth of dimension $k+2$. Moreover  the
normal bundle has the decomposition $\nc\cong {\Cal O}_C(k)
 \oplus {\Cal O}_C(k)$.
 \endproclaim

As another consequence of Theorem 2.3  we will show that  contact curves in
$\cpc$ can not be
complete intersections.

\proclaim{Corollary 2.6} Let $C$ be a  smooth contact curve in $\cpc$. If
$C$ is degenerate, then it is a straight line. If $C$ is non-degenerate, then
 it can not be a complete intersection.

\endproclaim

\demo{Proof}  In any case, by Lemma 2.1, we have the following short exact
sequence:
$$0\lr \omega_C(2)\lr \nc\lr {\Cal O}_C(2)\lr 0\tag 2.9$$
By Theorem 2.3, the above sequence does not splits.

 If $C$ is degenerate, then $C$ is contained in some hyperplane $\cpb\subset
  \cpc$. Let $d$ be the degree of $C$. On the one hand, the normal
  bundle has the decomposition \linebreak $\nc\cong {\Cal O}_C(1)\oplus {\Cal
O}_C(d)$.
  On the other hand, since  $\omega_C\cong {\Cal O}_C(d-3)$,  (2.9) becomes:
 $$0\lr {\Cal O}_C(d-1)\lr \nc\lr {\Cal O}_C(2)\lr 0\tag 2.10$$
 Therefore there is a non-trivial homomorphism $v: {\Cal O}_{C}(d)\lr
 {\Cal O}_C(2)$. Hence $d\leq 2$. If $d=2$, then $v$ is an isomorphism.
  Therefore (2.10) has to split. This is absurd since it does not split. This
   shows that $d=1$, i.e., $C$ is a straight line.

 If $C$ is non-degenerate, we will show that it is not a complete
 intersection. Suppose $C$ is a complete intersection, we write $C=S_n\cap
S_m$, where
$S_n$, resp. $S_m$ is surface of degree $n$, resp. $m$.   Since $C$ is
non-degenerate,
 both $n$ and $m$ are at least two. Suppose that $n\geq m$.
 Now  sequence (2.9) and the fact that
 $\nc\cong {\Cal O}_C(n)\oplus {\Cal O}_C(m)$ imply that
there is a non-trivial ${\Cal O}_C$-module homomorphism
$u: {\Cal O}_C(m) \lr
 {\Cal O}_C(2)$. Hence $m\leq 2$. Therefore $m=2$.  This implies that
$u$ is in fact an isomorphism. Hence (2.9) has to split. This is a
contradiction. Therefore $C$  is not a complete intersection. We are done.

\qed
\enddemo

\heading{\bf $\S$3. Contact Curves in $\cpc$}\endheading

   In this section we will  study contact curves in $\cpc$.  We will
 give a Pl\"ucker type formula for contact curves in $\cpc$.
This formula was predicted in [Br1] by Bryant. At the end of the section,
we will relate moduli space of contact line in $\cpc$ to
 the set of jumping lines for a null-correlation bundle.   In this section we
 will consider also {\it singular} contact curves. A singular curve $C$  is
  {\it contact} if local contact
forms vanish on $C_{reg}$, the smooth part of $C$.

     Since all contact structures on $\cpc$  are equivalent under
automorphisms of $\cpc$,  we will consider only  the distinguished contact
structure obtained from the twistor map $\cpc \overset T\to\lr S^4$.
 In  the affine coordinates $[1,z_1,z_2,z_3]$ of $\cpc$, the
local contact form is $\theta=d\,z_1-z_3\,d\,z_2+z_2\,d\,z_3$.  Bryant [Br1]
 provided a way to construct all contact curves in $\cpc$. He showed that
all contact curves in $\cpc$ are ``lifts" of curves in $\cpb$.

    First  of all,  let us make clear what the ``lift" means.  For a reduced
and irreducible curve $D\subset \cpb$, let
$D_{reg}$ be the smooth part of $D$. Then  points in $D_{reg}$ together
 with their tangent directions form a curve in $\ppb$. This gives a lift
of $D_{reg}$ to $\ppb$. Now the Zariski closure of this lift in $\ppb$
is the lift of $D$ to $\ppb$. Sometimes we call this the horizontal
lift of $D$, and we denote it  by $\tilde D$. It is clear that
  $\tilde D$ is smooth if $D$ has only unramified or simple cuspidal
  singularities.

 Bryant [Br1] defined a birational map $f: \ppb\lr \cpc$ sending
$\left(x,y,[\ll,\lll]\right)$ to $\big\lbrack\ll,x\ll-\dfrac{1}{2}y\lll,y\ll,
\dfrac{1}{2}\lll\big\rbrack$,
 where $(x,y)$ are coordinates on $\C^2\subset \cpb$ and
$[\ll,\lll]$ are fiber coordinates.

  Now we can state the theorem of Bryant [Br1].

\proclaim{Theorem 3.1}[Bryant] Let $C$ be a contact curve in $\cpc$. Then
$C$ is either a straight line or of the form  $f\left(\tilde D\right)$, where
$\tilde D\subset \cpb$ is the horizontal lift of a  reduced and
 irreducible plane curve $D$  of degree at least two.
\endproclaim

    Therefore to study contact curves in $\cpc$ is the same as to study
curves in $\cpb$ and their lifts to $\ppb$.  As we see in the
 first section, there is a natural contact structure on $\ppb$. The above
 birational map $f$ is {\it contact} in the sense that it
maps contact curves to contact curves.

 Next we offer a more homogeneous way of defining Bryant's map.
 We can think of $\ppb$ as the universal family of lines in $\cpb$, or as
  the flag manifold $F_{1,2}$. Therefore $\ppb$ sits naturally
  inside $\pp$, where $\cpbd$ is the dual of $\cpb$. In fact $\ppb$ is the set
of pairs $(p,\ell)$ such that
   $p\in \ell$, where $p$ is a point in $\cpb$ and $\ell$ is a line (which
is thought of as a point in $\cpbd$). Let $[x_i], \, [y_i]$ ($i=0,1,2$) be
 projective coordinates in $\cpb$ and $\cpbd$. Then $\ppb\subset\pp$ is defined
by the bi-homogeneous polynomial:
$$x_0y_0+x_1y_1+x_2y_2=0$$

Let $U_2=\{x_2\neq 0\}\subset \cpb$ be an affine open set
with affine coordinate $(x,y)=\left(x_0/x_2,x_1.x_2\right)$.   On $U_2$, we can
identify  the fiber coordinate $[\ll,\lll]$ with $[-y_0,y_1]$. Therefore
Bryant's map can be redefined as:
$$f\left([x_0,x_1,x_2], [y_0,y_1,y_3]\right)=
[2x_2y_0,2x_0y_0+x_1y_1,2x_1y_0,-x_2y_1] \tag 3.1$$
Consider two points $p_{\infty}=[1,0,0]\in \cpb$ and $p^*_{\infty}=[0,0,1]\in
\cpbd$. Then
$(p_{\infty}, p_{\infty}^*)$ is a point in $\ppb$. Let $\li\subset \cpb$ be the
line  dual to $p_{\infty}^*\in \cpbd$, and
 $\lid\subset \cpbd$ be the line dual to $p_{\infty}\in \cpb$.
Then it is clear that the $f$ is not defined precisely
along the curve $\li\cup \lid$.

 Let $p_1:\ppb\lr \cpb$ and $p_2:\ppb\lr \cpbd$ be  the two projections.
Note that if $D\subset \cpb$ is a
reduced and irreducible plane curve and $\tilde D$ is its lift, then
$p_1\left(\tilde D\right)=D$, and $p_2\left(\tilde D\right)=D^*$, where
$D^*$ is the dual curve of $D$. We say that a plane curve $D$ is {\it good}
 if $p_{\infty}\notin D$ and $p^*_{\infty}\notin D^*$. The second
condition is equivalent to the fact the $D$ doest not tangent to the line
$\li$.  Also $D$ is good if and only if $D\cap\left(\li\cup\lid\right)
=\emptyset$.  If $D$ is good and unramified, then  $\tilde D$ is smooth, hence
$f\left(\tilde D\right)$ is also smooth. Note that $D^*$ may be a point. This
happens exactly when $D$ is a line in $\cpb$. To make the following theorem
hold for lines also, we understand the degree of $D^*$  as zero if $D^*$ is a
point.  Now we get
the following theorem immediately.

\proclaim{Theorem 3.2} Let $C\subset \cpc$ be a contact curve of
degree $d$ and  geometric genus $g$, which is  obtained
from a  {\it good} plane curve $D$ of degree $n$.
 Then
$$\aligned
d &= n+n^*\\
g &= g(D)
\endaligned\tag 3.1$$
where $n^*$ is the degree of the dual curve of $D$ and $g(D)$ is the
geometric genus of $D$.
\endproclaim
\demo{Proof}  Let $U=\ppb\setminus\left(\ell_{\infty}\cup
\ell_{\infty}^*\right)$.
  By the definition of the Bryant's map, it is clear that
   $f^*{\Cal O}_{\cpc}(1)\vert_U\cong p_1^*{\Cal O}_{\cpb}(1)\otimes p_2^*{\Cal
O}_
  {\cpbd}(1)\vert_U$.  Since $D$ is {\it good}, $\tilde D$ is contained in $U$.
  Therefore
 the degree of $C\subset \cpc$ is the same as
 $p_1^*{\Cal O}_{\cpb}(1)\cdot \tilde D+p_2^*{\Cal O}_
  {\cpbd}(1)\cdot \tilde D=n+n^*$, where $n^*$ is the degree of the dual
  of $D$. Since $C$ and $D$ are obviously isomorphic to each other,
  $g(C)=g(D)$.
  Hence we are done.

\qed
\enddemo

In particular, if $D$ has only traditional singularities (see
[GH] for definitions), then the
following corollary is an immediate consequence of the above theorem
and the classical Pl\"ucker formula.
\proclaim{Corollary 3.3} Let $C$ be a  contact curve in $\cpc$  obtained
from a {\it good} plane curve $D$ of degree $n$ with traditional singularities.
 Then

$$\aligned
d &= n^2-2\delta-3\kappa\\
g &=\dfrac{1}{2}(n-1(n-2)-\delta-\kappa
\endaligned \tag 3.2
$$
where
 $\kappa$, resp. $\delta$ is the number of cusps, resp. nodes of $D$.
\endproclaim

\noindent{\bf Remark:} The requirement that $D$ is {\it good} simplifies
 the situation a lot, but it is not essential. By studying the resolution
of the map $f$ carefully, we should be able to handle the case of non-good
curves.
In that case, the degree of $C$ will be less than $n+n^*$ due to the
 the singularities of the map $f$.

Before we go on,  let us recall that
 a {\it complex involution} on a complex manifold $Y$
  is  an automorphism $\sigma$ of $Y$ such that $\sigma^2=id.$.

 Let $\hc(d,g)$ be the
 space  of contact curves (even including singular ones) in $\cpc$ of degree
$d$ and  geometric genus $g$.
 The above theorem shows that
 a plane curve and its dual give two contact curves with the same
degree and genus. In this way, we get a natural complex involution on
 $\hc(d,g)$. That is:
\proclaim{Corollary 3.4} There is  a  natural complex involution
on $\hc(d,g)$ . \endproclaim
\demo{Proof}  There is a natural involution (denoted by $\sigma_0$) on $\ppb$
sending
$[x_0,x_1,x_2]$ to $[y_0,y_1,y_2]$. Let $\tilde D$ be the lift of $D$.
 Then it is clear that $\sigma_0\left(D\right)=\tilde {D^*}$, the
lift of $D^*$ from $\cpbd$. Let $C^*=f\left(\tilde {D^*}\right)$. Then $C^*$
is also a contact curve $\cpc$. We thus define $\sigma\left(C\right)=C^*$.
Then $\sigma$ is an  complex involution. Hence we are done.

\qed
\enddemo

We next  present an  example, in which
the contact geometry is related to
the  moduli space of instanton bundles on $\cpc$. We will study
the moduli space $\hc(1,0)$,  the space of contact lines in $\cpc$.

\noindent {\bf Example 4. }\hskip 5pt
 Recall that the moduli space of contact structures
on $\cpc$ is isomorphic to ${\Cal M}(0,1)$, the moduli space of
rank-two stable bundles with $c_1=0$ and $c_2=1$ (see  Example 2.).
 As it is shown in
[Ha] that ${\Cal M}(0,1)$ is isomorphic to the space of all
non-singular anti-symmetric $4\times 4$ complex matrices. Hence it
 can be identified with $\C\P^5-G(1,3)$. Note that $G(1,3)$ is  the space of
all straight lines in $\cpc$. It is isomorphic to $\Q^4$, the  smooth
hyperquadric in $\C\P^5$.  As noted in [Ha]  $\Q^4$ induces a duality between
 non-singular hyperplane sections of $\Q^4$ and points in
$\C\P^5-G(1,3)$. Given a point $\v\in \C\P^5-G(1,3)$, draw all the lines
 through $\v$ and tangent to $\Q^4$. The points where these lines tangent
 to $\Q^4$ actually lie on an unique  linear subspace $\C\P^4\subset \C\P^5$.
 This linear subspace gives a hyperplane section $H_{\v}\subset \Q^4$. Barth
 [Ba] showed that
$H_{\v}$ can be identified with the jumping set $Z_{\v}$ for
$N_{\v}$,  the stable bundle corresponds to $\v$ in  Example 2. However
 we will  show that $Z_{\v}$ can be identified with $\hv\subset G(1,3)$, the
space of contact lines with respect to $\v$.

\proclaim{Proposition 3.5} The set of contact lines $\hv$  can be identified
with
the set of jumping lines of  $Z_{\v}$. In particular, $\hv$ is isomorphic
to  a nonsingular quadrics $\Q^3\subset \C\P^4$, hence  it has a
complex involution and a real structure.
\endproclaim
\demo{Proof} Let  $C$  is line and  is contact with respect to  the given
contact
structure $\v$.  Then the sequence (1.7)  for $\cpc$  implies that we have the
following short exact sequence:
$$ 0\lr {\Cal O}_C(2) \lr L^{\bot}\vert_C \lr {\Cal O}_C\lr 0 \tag 3.3$$
where $L={\Cal O}_{\cpc}(-2)$ is the contact bundle
 of  $\cpc$, and $L^{\bot}$ is the bundle of vector fields perpendicular to
  local sections of $L$. Since $H^1\left({\Cal O}_C(2)\right)=0$, the
  above sequence splits. Hence $L^{\bot}\vert_C\cong
  {\Cal O}_C\oplus {\Cal O}_C(2)$.  The  definition of $N_{\v}$ implies that
  $N_{\v}=L^{\bot}\vert_C\otimes {\Cal O}_{\cpc}(-1)$. Hence $N_{\v}\cong
  {\Cal O}_C(-1)\oplus {\Cal O}_C(1)$. Therefore $C$ is a jumping
  line for $N_{\v}$. So $\hv\subset Z_{\v}$. However $\text{dim}\,
\hv=h^0\left({\Cal O}_C(2)\right)=3$, which is the
dimension of $Z_{\v}$. Since $Z_{\v}$ is smooth and irreducible, we get
 $\hv=Z_{\v}$. Hence we are done.

   \qed
\enddemo

  We close this section by posing a question. Fix a pair of integers $(d,g)$
such that $d\geq g+3$.
  Ein [Ei] showed that the  Hilbert scheme $\hb(d,g)$ of smooth curves
in $\cpc$ with degree $d$ and genus $g$ is irreducible. Given
a contact structure $\v$ on $\cpc$, we have a closed subscheme $\hv(d,g)$ of
$\hb(d,g)$ consisting of  {\it smooth}  contact curves. Then it is natural to
ask:

\proclaim{Question 1} If   $\hv(d,g)$ is non-empty,  is  it
irreducible?.\endproclaim

By Proposition 2.2, $\hv(d,g)$ is smooth. Then if it is connected then it is
irreducible. Hence the above question is equivalent to asking if
$\hv(d,g)$ is connected. We remarked in the previous section  that $\hv(d,g)$
is a vanishing locus
 of a section $s_{\v}$ of a vector bundle $\E_{\v}$ over $\hb(d,g)$. One may
hope to
use Fulton-Lazarsfeld's connectivity theorem [FL]. But there are two
esential difficulties that make this approach impossible. First of
all, $\hb(d,g)$ is not projective. Secondly, most importantly, the vector
bundle
 $\E_{\v}$ is {\it not} ample. We have to take some other approach.

\noindent{\bf Remarks on higher-dimensional contact manifolds}

     Suppose that $n\geq 2$ and $X$ is a $(2n+1)$-dimensional  complex contact
manifold with a
	 contact line bundle $L$. Then in this situation, we can study
	  contact $n$-dimensional submanifolds in $X$.
	  Some of the general results in this  paper can be generalized  to
	  this situation without suitable modifications.
	  For example, Theorem 1.2, 1.4
	  and 2.3 can be generalized properly. However, Theorem 1.6 and
	  Theorem 1.8
	  can not. This is because  in higher dimensions we don't have a good
understanding
	  of extremal rays and pluri-canonical systems as we do in
	  three-dimensional  case. For example, we don't even  know if a
	  projective manifold $X$ having an extremal ray with length $\text{dim}
	  (X)+1$ is isomorphic to a  complex projective space.  In fact, this
	  is  an open conjecture in Mori's theory. However if $X$ is Fano, then
	   $X$ is isomorphic to the projectivized tangent bundle of
	   a projective space provided that Picard number of $X$
	   is at least two.  This is not difficulty to  prove
	   using claims of Wisniewski [Wi] (at the end of
	   the paper) and techniques of this paper.  We would like to ask the
	   following question:

	   \proclaim {Question 2} Let $X$ be a  contact $(2n+1)$-dimensional Fano
manifold
	   with Picard number $\rho(X)=1$. Is it true that  $X\cong \C\P^{2n+1}$?.
	   \endproclaim

\enddocument